\documentclass[preprint,3p]{elsarticle}

\usepackage{epstopdf}
\usepackage{pdfpages}
\usepackage{amssymb}
\usepackage{empheq}
\usepackage{cases}
\usepackage{amsthm,amsmath}
\usepackage{caption,lipsum}
\usepackage{stmaryrd}
\usepackage{graphicx,epsfig}
\usepackage{tabularx}
\usepackage{adjustbox}
\usepackage{color}
\usepackage{hyperref}
\usepackage{empheq,float}
\usepackage{booktabs,makecell,multirow,tabularx}
\usepackage{hhline}
\DeclareGraphicsExtensions{-eps-converted-to.pdf}
\newtheorem{theorem}{Theorem}[section]

\theoremstyle{definition}

\theoremstyle{remark}
\newtheorem{remark}[theorem]{Remark}
\newtheorem{thm}{Theorem}[section]

\newtheorem{lem}[thm]{Lemma}
\newtheorem{prop}[thm]{Proposition}
\newtheorem{defn}[thm]{ \bf{Definition}}
\newtheorem{conj}[thm]{ \bf{Conjecture}}

\newcommand{\EQ}[1]{\begin{align*}\begin{split} #1 \end{split}\end{align*}}
\newcommand{\EQn}[1]{\begin{align}\begin{split} #1 \end{split}\end{align}}

\newcommand{\Del}[1]{}

\def\pd{\partial}

\newcommand{\N}{{\mathbb N}}

\newcommand{\R}{{\mathbb R}}
\newcommand{\C}{{\mathbb C}}

\newcommand{\fe}{{\mathrm{e}}}

\def\De{\Delta}

\numberwithin{equation}{section}

\journal{}

\begin{document}
\begin{frontmatter}

\title{On blowup solution in NLS equation under dispersion or nonlinearity management}



\author[Chs]{Jing Li}\ead{lijingnew@csust.edu.cn}
\address[Chs]{School of Mathematical and Statistics, Changsha University of Science and Technology, 410114 Changsha, China}
\author[Gzh]{Cui Ning}\ead{cuining@gduf.edu.cn}
\address[Gzh]{School of Financial Mathematics and Statistics, Guangdong University of Finance, 510521 Guangzhou, China}
\address[Whu]{School of Mathematics and Statistics, Wuhan University, 430072 Wuhan, China}
\author[Whu,Hubei]{Xiaofei Zhao\corref{5}}\ead{matzhxf@whu.edu.cn}
\address[Hubei]{Computational Sciences Hubei Key Laboratory, Wuhan University, 430072 Wuhan, China}\cortext[5]{Corresponding author, URL: http://jszy.whu.edu.cn/zhaoxiaofei/en/index.htm.}

\begin{abstract}
In this paper, we study the dispersion-managed nonlinear Schr\"odinger (DM-NLS) equation
$$
i\partial_t u(t,x)+\gamma(t)\Delta u(t,x)=|u(t,x)|^{\frac4d}u(t,x),\quad x\in\R^d,
$$
and the nonlinearity-managed NLS (NM-NLS) equation:
$$
i\partial_t u(t,x)+\Delta u(t,x)=\gamma(t)|u(t,x)|^{\frac4d}u(t,x), \quad x\in\R^d,
$$
where $\gamma(t)$ is a periodic function which is equal to $-1$ when $t\in (0,1]$ and is equal to $1$ when $t\in (1,2]$.
The two models share the feature that the focusing and defocusing effects convert periodically.
For the  classical focusing NLS, it is known that the initial data
$$
u_0(x)=T^{-\frac{d}{2}}\fe^{i\frac{|x|^2}{4T}
-i\frac{\omega^2}{T}}Q_\omega\left(\frac{x}{T}\right)
$$
leads to a blowup solution
$$
(T-t)^{-\frac{d}{2}}\fe^{i\frac{|x|^2}{4(T-t)}
-i\frac{\omega^2}{T-t}}Q_\omega\left(\frac{x}{T-t}\right),
$$
so when $T\leq1$, this is also a blowup solution for DM-NLS and NM-NLS which blows up in the first focusing layer.

For DM-NLS, we prove that when $T>1$, the initial data $u_0$ above does not lead to a finite-time blowup and the corresponding solution is globally well-posed. For NM-NLS, we prove the global well-posedness for $T\in(1,2)$ and we construct solution that can blow up at any focusing layer. The theoretical studies are complemented by extensive  numerical explorations towards understanding
the stabilization effects in the two models and addressing their difference.
\end{abstract}


\end{frontmatter}

\section{Introduction}
\vskip .5cm

\subsection{Background}
The so-called dispersion management is an effective technique in fiber-optical communications. It can be implemented by several means in practice, e.g., the periodic concatenation of two different species of optical fibers with opposite values of the group-velocity dispersion \cite{GVD1,GVD2} or the insertion of
periodic nonlinear effects with alternating focusing and defocusing layers \cite{NM1,NM2}.
Mathematically, the basic model that these studies amount to is a $d$-dimensional ($d\in\N_+$) dispersion-managed nonlinear Schr\"odinger (DM-NLS) equation \cite{Saut,Malomed-book,Murphy}:
\EQn{\label{DM-NLS model}
i\partial_t u(t,x)+\gamma(t)\Delta u(t,x)=|u(t,x)|^{p-1}u(t,x),\quad
 x\in \mathbb{R}^d,\ t>0,
}
where $u=u(t,x):(0,\infty)\times\R^d\to\C$ is the unknown, $p>1$ is a given parameter, and the dispersion-map $\gamma(t)=\gamma(t+t_0)$ is a given periodic and piecewise-constant function in $t$:
\EQn{
	\label{DM}
\gamma(t):=\left\{
\aligned
-\gamma_-,\quad\, 0< t\leq t_*,\\
\gamma_+,\quad t_*< t\leq t_0,
\endaligned
\right.
	}
with some constants $\gamma_\pm>0$ and $t_0>t_*>0$. Such type of coefficient was introduced originally in fiber optics to control the accumulated dispersion in a traveling solitary wave, and so keep the wave narrowed in the focusing layer. This enhances the stability of dynamics in the physical system, which was found  particularly  important and useful in the applications like  signal processing, e.g., the suppression of Gordon-Haus jitter \cite{GVD2}.

\begin{figure}[h!]
  \centering
  \includegraphics[width=8.5cm]{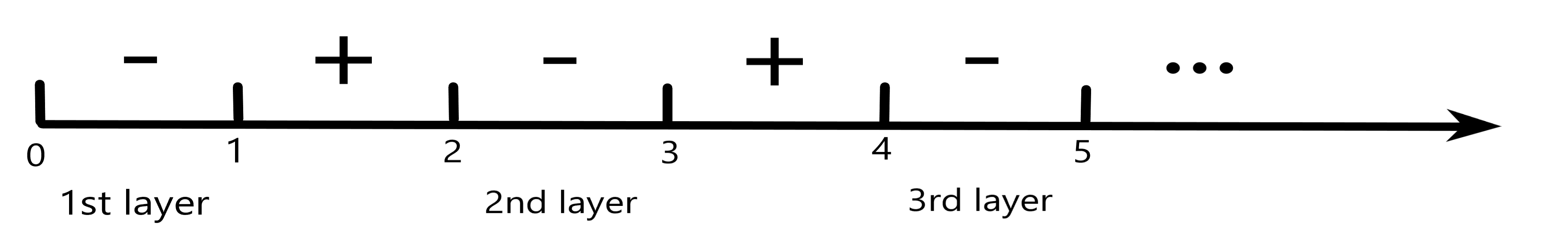}
  \caption{Illustrative example of the focusing layers.}\label{fig:layer}
\end{figure}

Apparently,  the DM-NLS (\ref{DM}) can be interpreted as a consecutively composition of the focusing NLS and the defocusing NLS in time, and they are referred as the layers. With $t_0=t_*=1$ in (\ref{DM}), a conceptual illustration is shown in Figure \ref{fig:layer}.
The model (\ref{DM})  conserves the mass for all times, i.e.,
$$\|u(t,\cdot)\|_{L^2}\equiv\|u(0,\cdot)\|_{L^2},\quad t\geq0,$$
and conserves the energy:
\begin{equation}\label{energy}
E_\gamma(t):=\int_{\mathbb{R}^d}\left[\frac{1}{2}|\nabla u|^2+\frac{\gamma(t)}{p+1}|u(t,x)|^{p+1}\right]dx
\end{equation}
piecewise in time, i.e., $\forall n\in\N$,
\begin{equation*}
  E_\gamma(s+nt_0)\equiv\left\{\begin{split}
                                  E_\gamma^+(nt_0),\qquad & 0< s\leq t_*, \\
                                   E_\gamma^+(t_*+nt_0),\qquad & t_*< s\leq t_0,
                               \end{split}\right.
\end{equation*}
with $E_\gamma^+(t)=\lim_{s\to t^+}E_\gamma(s)$ denoting the right limit. At the interface time: $t=nt_0$ or $t=t_*+nt_0$, the energy will exhibit a sudden change which makes it overall a step function, and the solution of (\ref{DM-NLS model}) will be discontinuous in time at the interface.

The DM-NLS model (\ref{DM-NLS model}) has attracted various investigations from aspects of physics, mathematics and numerics, e.g., \cite{2d-dm,Biswas-Book,dm-numer1,dm-numer2,Malomed-book} and the references therein. Although the majority physical studies concern the model (\ref{DM-NLS model}) in one space dimension with cubic nonlinearity, mathematically we are interested in the dramatic effect brought by the dispersion management to the NLS model in general.
Along the mathematical side, the consecutive switching between focusing and defousing layers
makes the dynamics in  (\ref{DM-NLS model}) very different from the classical NLS models \cite{Taobook}.
The very first attempt to understand the well-posedness theory of the Cauchy problem (\ref{DM-NLS model}) was made by Antonelli, Saut and Sparber in \cite{Saut}. They established the local well-posedness in $H^1$ for the energy-subcritical case and the
global well-posedness in $L^2$ for the mass-subcritical case.
The well-posedness theory in $H^1$ for the mass-critical case was established under a sharp criterion with the existence of  possible finite-time blowup. Also in their work, the fast dispersion management limit, i.e., $t_0\to0$ in (\ref{DM}), was investigated where an effective average model was obtained.  A recent work of Murphy and Hoose \cite{Murphy} established the global-in-time Strichartz estimates for
the linear case of the DM-NLS model. For the cubic DM-NLS case in three-dimension, they established a small-data scattering result and constructed a blowup solution in the first focusing layer.

A twin model imposes the dispersion-map (\ref{DM}) alternatively in front of the nonlinearity which introduces the nonlinearity-managed NLS equation (NM-NLS) \cite{Berge,Malomed-book,Saito,NM2}:
\EQn{\label{NM-NLS model}
i\partial_t u(t,x)+\Delta u(t,x)=\gamma(t)|u(t,x)|^{p-1}u(t,x),\quad
 x\in \mathbb{R}^d,\ t>0.
}
Such model shares quite some similarities as the DM-NLS: it also consists of focusing and defocuing layers as in Figure \ref{fig:layer}; the mass and energy of (\ref{NM-NLS model}) are preserved in the same way as in DM-NLS. However, both the mathematical and physical mechanism of the two models are different.
The concept of nonlinearity management originates from \cite{NM1}, which has also been used to stabilize the light signal transmissions in optical media \cite{Berge,NM stabilize,NM2}. Moreover, the NM-NLS (\ref{NM-NLS model}) is particularly relevant  for studying the Feshbach resonance \cite{NMadd0,Feshbach0,Feshbach1,VMPG1,Saito} in Bose-Einstein condensates, where the change of sign in $\gamma(t)$ can be experimentally realized by a magnetic field. Along the mathematical aspect of (\ref{NM-NLS model}), we emphasize that it does not have any equivalence with (\ref{DM-NLS model}) and it lacks of rigorous analysis. So far,  we are only aware of the works \cite{NM-theory1,NM-theory2} that analyzed (\ref{NM-NLS model}) under a continuous version of the dispersion map.

\subsection{Main results}
In this work, we shall analytically study the dynamics of both models, i.e., (\ref{DM-NLS model}) and (\ref{NM-NLS model}), complemented by some numerical explorations in the end.
The aim is to address the possible blowup phenomenon occurring in the second or latter focusing layer  (as depicted in Figure \ref{fig:layer}) of (\ref{DM-NLS model}) or (\ref{NM-NLS model}).

To simplify the notations and establish rigorous theorems, we restrict to analyze the initial value problem of the following `normalized'  mass-critical DM-NLS that will be shorted as $(NLS)_\gamma$:
\begin{equation}
	\label{eq:nls}
	\left\{ \aligned
	&i\pd_t u + \gamma(t)\De u = |u|^{\frac4d} u, \quad x\in\mathbb{R}^d,\ t>0,\\
	& u(0,x) = u_0(x),\quad x\in\mathbb{R}^d,
	\endaligned
	\right.
\end{equation}
where
\begin{equation*}
\gamma(t)=\left\{ \begin{split}
&-1, \quad\quad\quad 0<t\leq 1,\\
&+1,  \quad\quad\quad 1<t\leq 2,
\end{split}\right.
\end{equation*}
periodically extended in time with  period $t_0=2$.
The positive sign ``$\gamma(t)=1$'' in \eqref{eq:nls} denotes the defocusing source $(NLS)_+$ , and the negative sign ``$\gamma(t)=-1$" denotes the focusing one $(NLS)_-$.  The `normalized' NM-NLS model that will be shorted as $(NLS)^\gamma$ reads:
\begin{equation}
	\label{NLS2}
	\left\{ \aligned
	&i\pd_t u +\De u = \gamma(t)|u|^{\frac4d} u, \quad x\in\mathbb{R}^d,\ t>0,\\
	& u(0,x) = u_0(x),\quad x\in\mathbb{R}^d,
	\endaligned
	\right.
\end{equation}
with ``$\gamma(t)=1$" in \eqref{NLS2} denoting the defocusing source $(NLS)^+$  and ``$\gamma(t)=-1$" denoting the focusing one $(NLS)^-$.

The conservation laws of (\ref{eq:nls}) or (\ref{NLS2}) read the same as before with $p=\frac4d+1$ in (\ref{energy}).
Under the normalization, for both models, all the time intervals $(2n,2n+1)$ for any $n\in\N$ now correspond to the focusing layers. At the first (focusing) layer,  (\ref{eq:nls}) and (\ref{NLS2}) have just a sign difference, while at the first defocusing layer, the two equations become the same, which means that the two models are not equivalent under a simple complex conjugate nor are connected via any other transformations as far as we know.  We remark that starting with the focusing or defocusing section in the model will not make any essential differences in our studies.

To present our \textbf{main theoretical results}, let us further introduce the following notations and review some classical results for the standard NLS.

\begin{defn}[Symmetry group]
For some function $f(x,t)$ and constants $\omega\in \mathbb{R}^+, x_0\in\mathbb{R}, T\in\mathbb{R}$, we denote the scaling transform as
\begin{equation*}
\mathcal{T}_{\omega,x_0,\theta}f:=\fe^{i\theta}
\omega^{\frac{d}{2}}f\left(\omega(x-x_0),\omega^2t\right),
\end{equation*}
and we denote the pseudo-conformal transform as
\begin{equation*}
\mathcal{P}_Tf:=(T-t)^{-\frac{d}{2}}\fe^{i\frac{|x|^2}
{4(T-t)}}\bar{f}\left(\frac{x}{T-t},\frac{1}{T-t}\right),
\end{equation*}
with $\bar{f}$ the complex conjugate of $f$.
\end{defn}

It is known that both transforms preserve a defocusing/focusing NLS equation \cite{Taobook}, and so we have directly the following property.
\begin{prop}
For any $\omega\in \mathbb{R}^+, x_0\in\mathbb{R}, T\in\mathbb{R}$, if $u$ is a solution of $(NLS)_\gamma$ (or $(NLS)^\gamma$), then
$\mathcal{P}_T\mathcal{T}_{\omega,x_0,\theta}u$ is also a solution of $(NLS)_\gamma$ (or $(NLS)^\gamma$).
\end{prop}

For the focusing NLS equation:
\begin{align}\label{NLS-foc}
i\pd_t u - \De u = |u|^\frac{4}{d} u,
\end{align}
there exists a standing wave solution $R_\omega(t):=\fe^{i\omega^2t}Q_\omega$ for any $\omega>0$, where $Q_\omega=
Q_\omega(x):=\omega^\frac d2Q(\omega x)
$ with $Q(x)$ the ground state solution of
$$
-\De Q+Q-Q^{\frac4d+1}=0.
$$
Then  $\mathcal{P}_{T_0}R_{\omega_0}(t)$ is a solution of (\ref{NLS-foc}) for any $\omega_0,T_0>0$, and it is known \cite{Merle-BU-1993} that
$$
\mathcal{P}_{T_0}R_{\omega_0}(t) \mbox{ blows up at } t=T_0.
$$
For simplicity, we denote $R(t)=R_1(t)=\fe^{it}Q$.
Owning to this fact,   Antonelli, Saut and Sparber in \cite{Saut} (Murphy and Hoose later in \cite{Murphy} for the inter-critical case) proved that by choosing $T_0\in (0,1)$, there exists a blowup solution $\mathcal{P}_{T_0}R_{\omega_0}$ for the DM-NLS equation \eqref{eq:nls},  which blows up in the first focusing layer.

Mathematically,  it would be natural to ask that for $T_0>1$,  whether the solution of \eqref{eq:nls} with the initial data $\mathcal{P}_{T_0}R_{\omega_0}(0)$ still blows up or not. More precisely, the question is whether the defocusing layer(s) will delay or even suppress the blowup.  The following result provides a theoretical answer.

\begin{thm}[DM-NLS result]\label{thm:main1}
For any $1<T_0<\infty$ and $\omega_0>0$, let $u_0=\mathcal{P}_{T_0}R_{\omega_0}(0)$, then the corresponding solution of the $(NLS)_\gamma$ \eqref{eq:nls}  is globally well-posed.
\end{thm}

The above theorem means that if $\mathcal{P}_{T_0}R_{\omega_0}(0)$ does not blow up in the first focusing layer, then it would never blow up in future. Thus, the dispersion management  could  indeed suppress the appearance of  blowups in the usual focusing NLS, which to some extent enhances the stability of the model and makes it more robust.  We will given the proof of Theorem \ref{thm:main1} in Section \ref{sec:proof1} and provide the numerical verification in Section \ref{sec:num}.

Let us move on to the NM-NLS model (\ref{NLS2}). Note that its focusing layer $(NLS)^-$ is a complex conjugate of (\ref{NLS-foc}), so $\overline{\mathcal{P}_{T_0}R_{\omega_0}(t)}$ is a blowup solution for $(NLS)^-$. On the one hand, we can show that the nonlinearity management in (\ref{NLS2}) could also help to suppress certain blowup, which is stated as the first part of result in the following theorem. Unfortunately,  due to the technical reason, we are only able to prove this for $T_0>1$ but not large. On the other hand, one interesting question that we are not able to address for DM-NLS is that whether there exists a solution that remains finite at the first few layers but blows up at a latter focusing layer. For the NM-NLS case, we can answer with a rigorous YES!
More precisely, we have the following statement given as the second part of the theorem, which implies that the solution of (\ref{eq:nls}) could blow up at an arbitrary focusing layer.



\begin{thm}[NM-NLS result]\label{thm:main2}
For the $(NLS)^\gamma$ (\ref{NLS2}),
the following two assertions hold:
\begin{itemize}
\item[(i)] For any $T_0\in(1,2)$, let $u_0=\overline{\mathcal{P}_{T_0}R_{\omega_0}(0)}$, then the corresponding solution of  \eqref{NLS2}  is globally well-posed.
\item[(ii)]
For any  $n\in \mathbb{Z}^+$ and any $T_*\in (2n, 2n+1)$, there exists a blowup solution $u_{T_*}$ of  \eqref{NLS2},  which blows up at $T_*$. 
\end{itemize}
\end{thm}

The above part (i) indicates the enhanced stability of the NM-NLS model. With such theoretically guaranteed stability, both DM-NLS and NM-NLS models may find perspective applications in broader areas.
The proof of Theorem \ref{thm:main2} will be given in Section \ref{sec:proof2} and verified in Section \ref{sec:num}.

We remark that the part (i) of Theorem \ref{thm:main2}  for $T_0\geq2$ as well as the counterpart of (ii) for the DM-NLS case may also hold, while our analytical tools are not enough to rigorously establish them.  They will be addressed by numerical explorations later. The difference between the obtained results in Theorem \ref{thm:main1} and Theorem \ref{thm:main2}, on the other hand, demonstrates the inherent  difference between the two NLS models.


\subsection{Some conjectures}
 Extensive numerical experiments will be done in Section \ref{sec:num} to complement the above theoretical results. That is to first verify the theorems and then explore general situations. Based on the numerical observations in Section \ref{sec:num} and many more from the physical literatures (e.g., \cite{VMPG1} and the references therein), we can draw  three formal conjectures that  further demonstrate  the stabilization effects in
the DM-NLS and NM-NLS models in the following. The purpose of the conjectures here is, on the one hand, summarize the numerical outcomes, and more importantly to draw mathematical attentions for future efforts of rigorous proof.

\begin{conj}\label{conj 1}
Let $u_0\in H^1(\R^d)$ such that the solution to  \eqref{NLS-foc} with $u(0,x)=u_0(x)$ blows up at some time $T>1$, then the solution of  \eqref{eq:nls} or \eqref{NLS2} with $u(0,x)=u_0(x)$ does not blow up in $[0,T]$. Thus, the blowups can  be at least delayed.
\end{conj}

\begin{defn}
Given an initial data $u_0$, we call the equations  \eqref{eq:nls} and \eqref{NLS2} dispersion-manageable, if their solutions under $u(0,x)=u_0(x)$ globally exist, and  there exist positive constants $c_0, C_0\in(0,\infty)$ such that
\begin{equation}\label{manageable}
c_0\le \sup\limits_{t\geq0} \big\|u(t,\cdot)\big\|_{L_x^\infty(\R^d)}\le C_0.
\end{equation}
\end{defn}

The dispersion-manageable property  implies that the solution is global but non-scattering. It behaves like a `solition'.
This is against the common dispersive  property of the Schr\"odinger equation.

\begin{conj}\label{Conjecture DM}
There exists a Schwartz function $u_0$ such that  $(NLS)_\gamma$ or $(NLS)^\gamma$  with $u(0,x)=u_0(x)$ is dispersion-manageable.
\end{conj}


\begin{defn}
We call an initial data $u_0$ is a non-scattering data if the solution to $(NLS)_-$ with $u(0,x)=u_0(x)$ does not
scatter.
\end{defn}

\begin{conj}\label{Conjecture DM 2}
Given a non-scattering data  $u_0\in H^1(\R^d)$,  then the following two assertions hold:
\begin{itemize}
\item[(i)] Fixing $t_0,t_1$, there exist $\gamma_\pm=\gamma_\pm(u_0)$;
\item[(ii)] Fixing $\gamma_+,\gamma_-$, there exist $t_0=t_0(u_0),t_1=t_1(u_0)$;
\end{itemize}
such that the solution to  \eqref{DM-NLS model} with $u(0)=u_0$ is dispersion-manageable.
\end{conj}

This conjecture implies that for any non-scattering data, suitably adjusting the system may lead to a ``soliton-like'' solution. 

The following are some further discussions.

(1) \textit{Energy behavior.} The energy $E_\gamma(t)$ is not generally conserved under the flow (it conserves in each separate layer). Then, the behavior of the energy may help us to understand the behavior of the solution. Numerical experiments given in Section \ref{sec:num} show that the energy $E_\gamma(t)$ grows in time even if the solution is global existence.  If the energy grows to infinity, then this in fact implies that the solution blows up at infinite time.

(2) \textit{Fast switching case.}  We may consider the equation in the rapidly varying case:
\begin{align}\label{eqs:nls-fs}
i\pd_t u_\varepsilon+ \gamma\left(\frac{t}{\varepsilon}\right)\De u_\varepsilon = |u_\varepsilon|^{\frac4d} u_\varepsilon.
\end{align}
That is,  we replace $\gamma(t)$ by $\gamma(\frac{t}{\varepsilon})$ in \eqref{eq:nls}  with $0<\varepsilon\ll 1$, then Theorems \ref{thm:main1} and \ref{thm:main2} still hold for any fixed $\varepsilon>0$.  As proved in \cite{Saut}, when $0<\varepsilon\ll 1$, the solution of \eqref{eqs:nls-fs} is close to an oscillatory (in time) function  which solves the equation
$$
i\pd_t u= |u|^{\frac4d} u.
$$
Then, the solution $u_\varepsilon$ may behave as a high-frequency oscillatory function when $\varepsilon$ is small.

(3) \textit{More general cases.} The proofs of Theorems \ref{thm:main1} and \ref{thm:main2}  rely on the classification lemma established by F. Merle in \cite{Merle-BU-1993}, which technically is only available for the mass-critical NLS under  the critical mass $\|Q\|_{L^2}$. For  more general cases, such as the mass-supercritical NLS or the mass-critical NLS with $E(u_0)<0$,  we expect similar results to be also valid.

\subsection{Organization}
The rest part of paper is organized as follows. In Section \ref{sec:lemma}, we present some useful lemmas that are key for the analysis. In Section \ref{sec:proof1} and Section \ref{sec:proof2}, we give the proofs of the two main theorems. Some numerical observations are made in Section \ref{sec:num}, and some concluding remarks are drawn in Section \ref{sec:con}.


\section{Some useful lemmas} \label{sec:lemma}
\vskip .5cm

Firstly, we need the following classification result from F. Merle \cite{Merle-BU-1993}, which plays a crucial role in our analysis.

\begin{lem}[F.Merle]\label{lem:classicification}
Suppose that $u_0\in H^1(\R^d)$ and $\|u_0\|_{L^2}=\|Q\|_{L^2}$. Let $u$ be the solution of $(NLS)_-$ with $u(0,x)=u_0(x)$.  Then  $u(t,x)$ blows up at the finite time $T_*$ iff there exist
$\omega, \theta, x_0$, such that
\begin{equation}
u=\mathcal{P}_{T_*}\mathcal{T}_{\omega, x_0, \theta}R.
\end{equation}
\end{lem}
\begin{remark}
According to Lemma \ref{lem:classicification}, we have the following blow-up solution for $(NLS)^{-}$ equation \eqref{NLS2},
\begin{align*}
u=\overline{\mathcal{P}_{T_*}\mathcal{T}_{\omega, x_0, \theta}R}.
\end{align*}
\end{remark}

Then, we define two auxiliary functionals:
$$
I(u):=\int_{\R^d}|x|^2|u|^2\,dx,
$$
and
$$
P(u):=\mbox{Im}\int_{\R^d}x\cdot\nabla u\,\overline{u}dx.
$$
To prove the main theorems, we also need the following two basic lemmas.

\begin{lem}\label{eqs:h} Let the function
\begin{equation*}
h_{T,\omega}(t)=\mathcal{P}_{T}R_\omega(t).
\end{equation*}
Then
\EQn{
	\label{eq:xh}
\aligned
I\big(h_{T,\omega}(t)\big)
&=\left(\frac{T-t}{\omega}\right)^2\big\|xQ\big\|_{L^2}^2,
\endaligned
	}
\end{lem}
and
\EQn{
	\label{eq:ph}
\aligned
P\big(h_{T,\omega}(t)\big)&=\frac{T-t }{2\omega^2}\|xQ\|_{L^2}^2.
\endaligned
	}
\begin{proof}
Recall the definition of $\mathcal{P}_{T}$, we can write
$$h_{T,\omega}(t)=\left(\frac{\omega}{T-t}\right)^{\frac{d}{2}}
\fe^{i\frac{|x|^2}{4(T-t)}-i\frac{\omega^2}{T-t}}
Q\left(\frac{\omega x}{T-t}\right).$$
For convenience, we denote $\frac{\omega x}{T-t}=x_{T,\omega}$.
Then, we have
\begin{align*}
I\big(h_{T,\omega}(t)\big)
&=\left(\frac{\omega}{T-t}\right)^d\>\int_{\R^d} \big|\,xQ(x_{T,\omega})\big|^2\,dx\\
&=\left(\frac{T-t}{\omega}\right)^2\>\int_{\R^d} \big|x_{T,\omega}\,Q(x_{T,\omega})\big|^2\,dx_{T,\omega}\\
&=\left(\frac{T-t}{\omega}\right)^2\big\|x\,Q\big\|_{L^2}^2.
\end{align*}
Furthermore,
\begin{align*}
\nabla h_{T,\omega}(t)=&\fe^{i\frac{|x|^2}{4(T-t)}
-i\frac{\omega^2}{T-t}}\>\left(\frac{\omega}{T-t}\right)^{\frac{d}{2}}
\left[\frac{ix}{2(T-t)}Q(x_{T,\omega})
+\frac{\omega}{T-t}\nabla Q(x_{T,\omega})\right].\\
\end{align*}
Thus,
\EQ{
\aligned
P(h_{T,\omega}(t))&=\mbox{Im}\int_{\R^d} x\cdot\nabla h_{T,\omega}(t)\,\overline{h_{T,\omega}(t)}dx\\
&=\left(\frac{\omega}{T-t}\right)^d\>\mbox{Im}\int_{\R^d} x\cdot
\left[\frac{ix}{2(T-t)}\,Q(x_{T,\omega})
+\frac{\omega}{T-t}\nabla Q(x_{T,\omega})\right]\,Q(x_{T,\omega})dx\\
&=\left(\frac{\omega}{T-t}\right)^d\>\mbox{Im}\int_{\R^d} x\cdot \frac{ix}{2(T-t)}Q(x_{T,\omega})\,Q(x_{T,\omega})\,dx\\
&=\frac{T-t}{2\omega^2}\int_{\R^d}|x_{T,\omega}|^2\,|Q(x_{T,\omega})|^2\,dx_{T,\omega}\\
&=\frac{T-t}{2\omega^2}\|xQ\|_{L^2}^2.
\endaligned
	}
This finishes the proof of the lemma.
\end{proof}

Note that the solution of (\ref{eq:nls}) or (\ref{NLS2}) can afford one time derivative, and by direct computing, we can have the following Virial identities.
\begin{lem}[Virial identities]\label{IP1}
The following identities hold.
\begin{itemize}
\item[(1)]
Let $u$ be the solution of $(NLS)_\gamma$, then
\begin{align*}
&\frac{d}{dt}I\big(u(t,\cdot)\big)=4\gamma(t)\mbox{Im}\int_{\R^d} x\cdot\nabla u\,\overline{u}dx=4\gamma(t)P\big(u(t,\cdot)\big);\\
&\frac{d}{dt}P\big(u(t,\cdot)\big)=2\gamma(t)\|\nabla u\|_{L^2}^2+\frac{2d}{d+2}\|u\|_{L^{\frac4d+2}}^{\frac4d+2}=4\gamma(t)E_{\gamma}\big(u(t,\cdot)\big).
\end{align*}
\item[(2)]
Let $u$ be the solution of  $(NLS)^{\gamma}$, then
\begin{align*}
&\frac{d}{dt}I\big(u(t,\cdot)\big)=4\mbox{Im}\int_{\R^d} x\cdot\nabla u\,\overline{u}dx=4P\big(u(t,\cdot)\big);\\
&\frac{d}{dt}P\big(u(t,\cdot)\big)=2\|\nabla u\|_{L^2}^2+\gamma(t)\frac{2d}{d+2}\|u\|_{L^{\frac4d+2}}^{\frac4d+2}=4E_{\gamma}\big(u(t,\cdot)\big).
\end{align*}
\end{itemize}
\end{lem}

The detailed calculations for the above are omitted for brevity.

\section{Proof of Theorem \ref{thm:main1}}\label{sec:proof1}
In this and the next section, we shall omit the spatial variable in the solution, i.e., $u(t)=u(t,x)$, for simplicity of notations.

We prove Theorem \ref{thm:main1} by contradiction. Assume that the solution blows up at some $T_*>1$. Since the solution obeys $(NLS)_+$ when $t\in [2j-1,2j]$ for $j\in \N^+$, it is impossible to blow up in this time region.
Hence, there exists some $k$ such that $T_*\in (2k,2k+1]$. Since the solution $u$ obeys the equation of  $(NLS)_-$ at $t\in [2k,2k+1]$, by applying Lemma  \ref{lem:classicification}, we have that there exist some $\omega>0, x_0\in \R, \theta\in\R$, such that
$$
u(t)=\mathcal{P}_{T_*}\mathcal{T}_{\omega, x_0, \theta}R (t)\quad  \mbox{with }\quad t\in [2k,T_*).
$$
Since $u_0$ is radial,  $u(t)$ is also radial for any $t\in [0,T_*)$. Thus, $x_0=0$. Moreover, the rotation has no influence on our analysis. Without loss of generality, we may set $\theta=0$. Hence, we may write
$$
u(t)=\mathcal{P}_{T_*}R_\omega (t)\quad  \mbox{with }\quad t\in [2k,T_*).
$$
From the definition of energy, we have
$$E_\gamma(t)=\frac12 \|\nabla u\|_{L^2}^2+\gamma(t)\frac{d}{2d+4}\|u\|_{L^{\frac4d+2}}^{\frac4d+2}.$$
Note that $\|u(t)\|_{L^2}=\|Q\|_{L^2}$  for both $\gamma(t)=-1$ and $\gamma(t)=1$, then by the Gagliardo-Nirenberg inequality, we always have
$$E_\gamma(t)\geq0.$$

In the following, we denote
$
P(t)=P\big(u(t)\big)
$
, $E_{-}=E_{-1}$ and $E_{+}=E_{+1}$
for short.

Since $u_0=\mathcal{P}_{T_0}R_{\omega_0}(0)=h_{T_0,\omega_0}$, we can conclude that
$$E_{-}(u_0)=E_{-}(h_{T_0,\omega_0})=\frac{1}{8\omega_0^2}\|xQ\|_{L^2}^2.$$
By Lemma \ref{IP1} and the conservation of energy, we have
\begin{align*}
 P'(t)=-4E_-(t)=-4E_{-}(u(2k))\leq0, \ t\in (2k,2k+1).
\end{align*}

If $u(t)$ blows up at $T_*$, we have $P(T_*)=0$. Because of the monotonicity of $P(t)$,
we have $P(2k+1)\leq0$.

On the other hand, from Lemma \ref{IP1}, it yields that
\begin{align*}
P(2k+1)
=P(2k)-4E_{-}(u(2k)),
\end{align*}
and
\begin{align*}
P(2k)
=P(2k-1)+4E_{+}(u(2k)).
\end{align*}
In particular, we have
\begin{align*}
P(1)=P(0)-4E_{-}(u_0).
\end{align*}
By Lemma \ref{eqs:h} and $T_0>0$, we have
\begin{align*}
P(1)=\frac{T_0}{2\omega_0^2}\|xQ\|_{L^2}^2-\frac{1}{8\omega_0^2}\|xQ\|_{L^2}^2=\frac{T_0-1}{2\omega_0^2}\|xQ\|_{L^2}^2>0.
\end{align*}
Then, we have
\begin{align*}
P(2k+1)=&P(2k-1)+4[E_{+}(u(2k))-E_{-}(u(2k))]\\
=&P(2k-1)+\frac{2d}{2d+4}\|u(2k)\|_{L^{\frac4d+2}}^{\frac4d+2}\\
=&P(1)+\sum_{j=1}^k\frac{2d}{2d+4}\|u(2j)\|_{L^{\frac4d+2}}^{\frac4d+2}>0.
\end{align*}

This is a contradiction, which implies that we have the global existence for the solution in $\R^+$. $\Box$

\section{Proof of Theorem \ref{thm:main2}}\label{sec:proof2}
\subsection{Proof Theorem \ref{thm:main2} (i)}
According to Lemma \ref{lem:classicification}, for some $\omega>0, x_0\in \R, \theta\in\R$, we have the corresponding blowup solution for $(NLS)^{-}$ equation \eqref{NLS2}
\begin{align*}
u(t)=\overline{\mathcal{P}_{T_*}\mathcal{T}_{\omega, x_0, \theta}R}.
\end{align*}

Without loss of generality, we may set $x_0=\theta=0$. Hence, we can write
\begin{align*}
u(t)=\overline{\mathcal{P}_{T_*}\mathcal{T}_{\omega}R}.
\end{align*}

Again, we prove it by contradiction.
Let the function
\begin{equation*}
H_{T,\omega}(t)=\overline{\mathcal{P}_{T}R_\omega(t)}.
\end{equation*}
By the direct computation, we get
\EQn{
	\label{eq:pH}
\aligned
P(H_{T,\omega}(t))=\mbox{Im}\int_{\R^d} x\cdot\nabla H_{T,\omega}(t)\,\overline{H_{T,\omega}(t)}dx&=\frac{t-T}{2\omega^2}\|xQ\|_{L^2}^2
\endaligned
	}
and
\EQn{
	\label{eq:xH}
\aligned
I(H_{T,\omega}(t))=\big\|xH_{T,\omega}(t)\big\|_{L^2}^2
&=\left(\frac{T-t}{\omega}\right)^2\big\|xQ\big\|_{L^2}^2.
\endaligned
	}
 Since $u_0=\overline{\mathcal{P}_{T_0}R_{\omega_0}(0)}$, then we have
$$u(1)=H_{T_0,\omega_0}(1).$$
By \eqref{eq:pH}, we have
\begin{align*}
P(1)=\frac{1-T_0}{2\omega_0^2}\|xQ\|_2^2.
\end{align*}
On the other hand, from Lemma \ref{IP1}, for any $t\in(1,2)$, we have
\begin{align*}
P'(t)&=4E_+(u(t))=4E_+(u(1))=2 \|\nabla u(1)\|_{L^2}^2+\frac{2d}{d+2}\|u(1)\|_{L^{\frac4d+2}}^{\frac4d+2}\\
&=\frac{1}{2\omega_0}\|xQ\|_{L^2}^2+4\left(\frac{\omega_0}{T_0-1}\right)^2\|\nabla Q\|_{L^2}^2.
\end{align*}
Then, for $T_0\in(1,2)$, we have
\begin{align*}
P(T_0)&=P(1)+(T_0-1)P'(t)\\
&=\frac{1-T_0}{2\omega_0^2}\|xQ\|_2^2+(T_0-1)\left[\frac{1}{2\omega_0}\|xQ\|_{L^2}^2+4\frac{\omega_0}{T_0-1}\|\nabla u\|_{L^2}^2\right]\\
&=4(T_0-1)\left(\frac{\omega_0}{T_0-1}\right)^2\|\nabla u\|_{L^2}^2>0.
\end{align*}
For any $t\in \mathbb{R}^+$, we know that $P'(t)=4E_\gamma(u(t))\geq 0$. Then, for any $T>T_0$, we obtain
$$
P(T)\geq P(T_0)>0.
$$
If $u(t)$ blows up at some $T>T_0$, we have $P(T)=0$. This is a contradiction. $\Box$

\subsection{Proof Theorem \ref{thm:main2} (ii)}\label{subsection thm2ii}
Let $T_*\in (2n,2n+1)$ for some $n\in \N^+$. We claim that there exists a solution $u$ satisfying
\begin{align}\label{sol-u}
u(t)=\overline{\mathcal{P}_{T_*}R_{\omega} (t)}\quad  \mbox{with }\quad t\in (2n,2n+1].
\end{align}
If so, then $u(t)$ blows up at $t=T_*$.

To prove the claim, it is sufficient to show that there exists an initial data $u_0$ for (\ref{NLS2}) such that the corresponding
solution $u$ of the NM-NLS obeys \eqref{sol-u}. This is equivalent to show that there exists a unique solution $\tilde u(t)$  on the time interval $[0,2n]$ for the following problem:
\begin{equation}
	\left\{ \aligned
	&i\pd_t {\tilde u} +\De \tilde u = \gamma(2n-t)| \tilde u|^{\frac4d} \tilde u, \quad x\in\mathbb{R}^d,\ t>0,\nonumber\\
	& \tilde u(0) =\overline{u(2n)}=\mathcal{P}_{T_*}R_{\omega} (2n),\quad x\in\mathbb{R}^d.
	\endaligned
	\right.
\end{equation}
Clearly, if it holds, then $u(t)=\overline{\tilde u(2n-t)}$.

To simplify the notations, we may set $x_0=\theta=0$ by spatial and rotation transformations.
In particular, we have  $$u(2n)=\overline{\mathcal{P}_{T_*-2n}R_{\omega^*}(0)}=H_{T_*,\omega_*}(2n).$$

%

Indeed, from the local theory, there exists $\delta>0$, such that $\tilde u(t)\in C((0,\delta),H^1)$. If there exists $T\in (0,2n)$, such that $\tilde u(t)$ blows up at $T$ with $\tilde u(0)=\overline{H_{T^*,\omega^*}(2n)}$, then we have
$$P(\tilde u(T))=0,$$
and
$$P(\tilde u(0))=P(\overline{H_{T_*,\omega_*}(2n)})=\frac{T_*-2n}{2\omega_*}\|xQ\|_{L^2}^2>0.$$
On the other hand, by Lemma \ref{IP1} and the Gagliardo-Nirenberg inequality, for $t\in(0,2n)$, we know that
$$P'(\tilde u(t))=2\|\nabla \tilde u\|_{L^2}^2+\gamma(2n-t)\frac{2d}{d+2}\|\tilde u\|_{L^{\frac4d+2}}^{\frac4d+2}\geq0.$$
Therefore, we find
$$P(\tilde u(T))>0,$$
which is a contradiction.

%
%

This finishes the proof of Theorem \ref{thm:main2} (ii). $\Box$

\begin{remark}
Here we briefly explain the reason why we cannot get the same assertion for DM- and NM-NLS models in Theorems \ref{thm:main1} \& \ref{thm:main2}.
The essential difference comes from the fact  that  $P(u(t))$ under the DM-NLS case  and  the NM-NLS case exhibits different behaviors as a function of time. Indeed from Lemma \ref{IP1} (1), $P'(u(t))$ for DM-NLS  changes the sign at each layer, while it does not for NM-NLS. This causes technical difficulties to construct blowup solutions for DM-NLS and  to prove the general global well-posdeness for NM-NLS.
\end{remark}

\section{Numerical observations}\label{sec:num}
\vskip .5cm

In this section, we shall perform some numerical simulations of the dynamics in the DM-NLS equation and in the NM-NLS equation. For discretizations of the two models, we would first truncate the whole space $\R^d$ to a finite interval with periodic boundary conditions and then apply the numerical schemes from \cite{Zhao} under fine mesh for accurate computations. 

In the literature, intensive numerical tests have already been carried out towards understanding the stabilization of waves in DM-NLS and NM-NLS. Many of them are performed in two-dimensional or three-dimensional cases, e.g., \cite{VMPG1} and the references therein, for controlling the blowups. Here, our presentations in this section mainly have two purposes. The first purpose is to support the proved theoretical results under  the precise setups in the theorems. The second purpose is to self-consistently support the drawn conjectures which cannot be proved mathematically at this moment. Indeed, some of our numerical results share the same findings from the literature. We highlight that our results particularly cover also the one-dimensional case of DM-NLS and NM-NLS, where both stabilization and blowup scenario are constructed, and the energy behaviour is included as well.

\begin{figure}[t!]
\centering
\psfig{figure=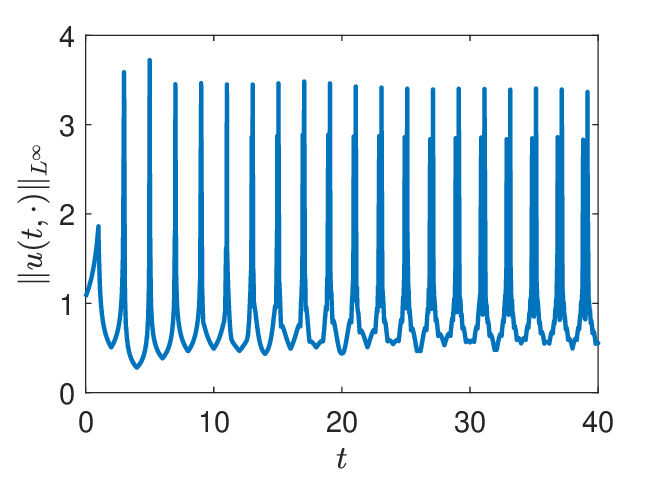,width=0.45\textwidth}
\psfig{figure=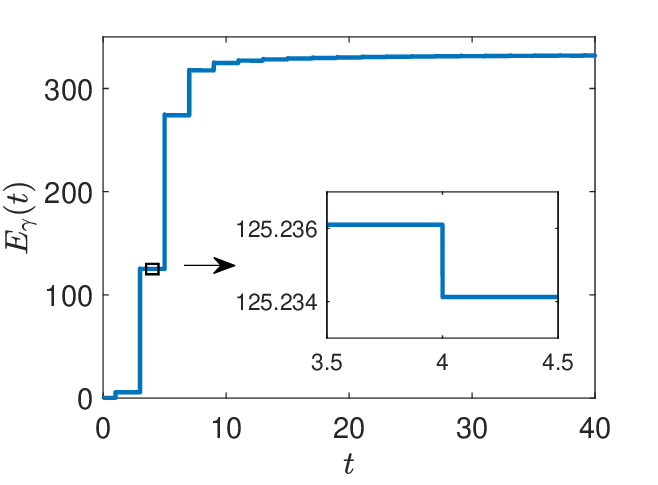,width=0.45\textwidth}
\caption{DM-NLS: evolution of the maximum value $\|u(t,\cdot)\|_{L^\infty}$ (left) and the energy $E_\gamma(t)$ (right).
}
\label{fig:1}
\end{figure}

\begin{figure}[t!]
\centering
\psfig{figure=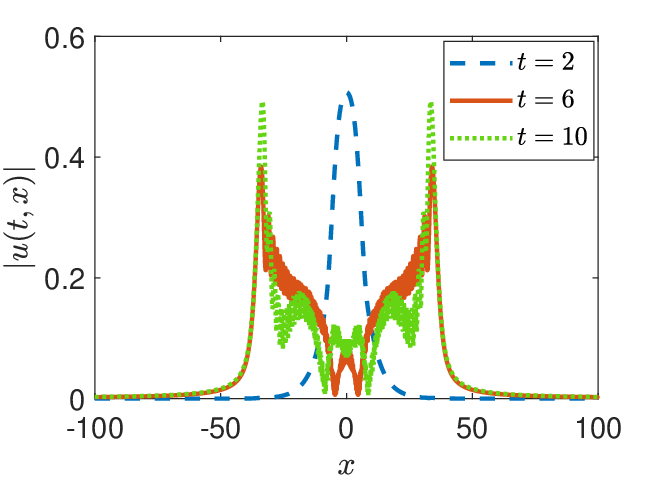,width=0.45\textwidth}\psfig{figure=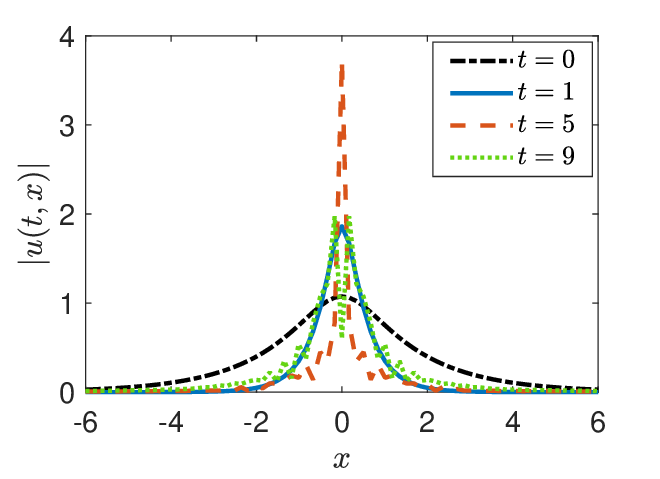,width=0.45\textwidth}\\
\psfig{figure=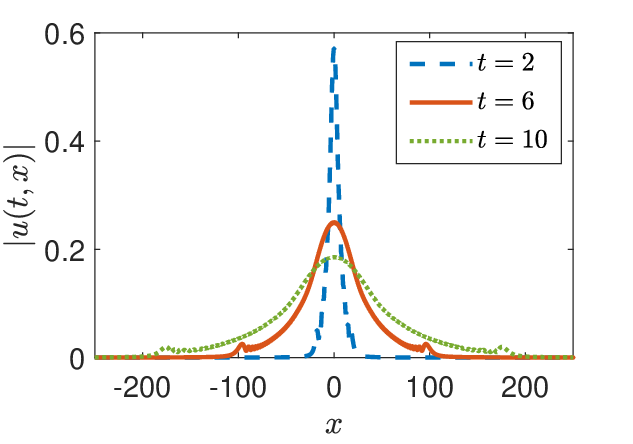,width=0.45\textwidth}\psfig{figure=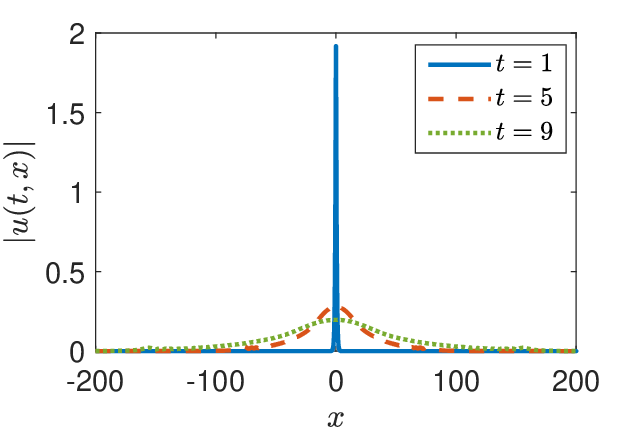,width=0.45\textwidth}
\caption{Profiles of the solution $|u(t,x)|$ at different time: DM-NLS (up) and NM-NLS (down).
}
\label{fig:1a}
\end{figure}

\begin{figure}[t!]
\centering
\psfig{figure=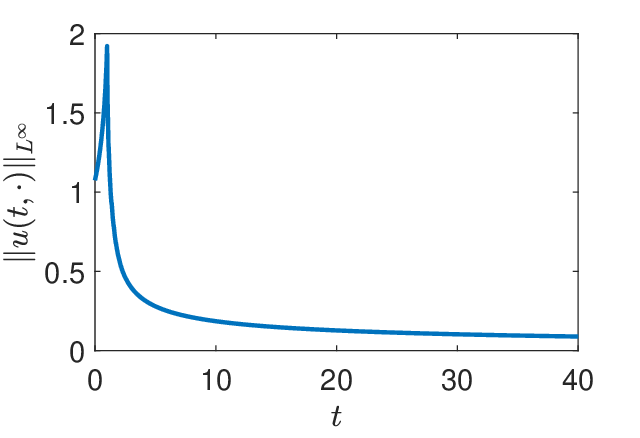,width=0.45\textwidth}
\psfig{figure=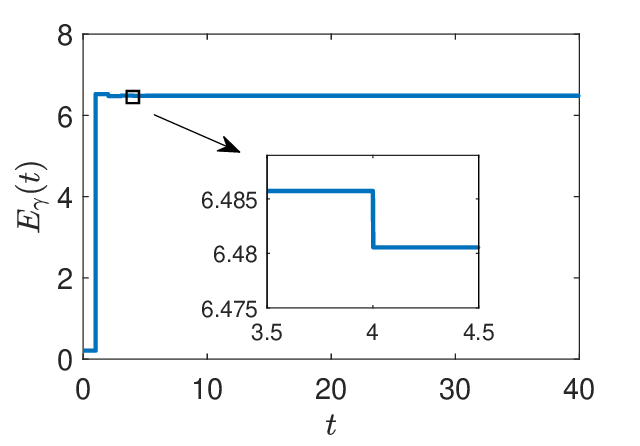,width=0.45\textwidth}
\caption{NM-NLS: evolution of the maximum value $\|u(t,\cdot)\|_{L^\infty}$ (left) and the energy $E_\gamma(t)$  (right).
}
\label{fig:2}
\end{figure}

\subsection{Verification of theoretical results}
Let us first verify the results  given in the two theorems.
For simplicity, we consider here the one-dimensional case, i.e., $d=1$ in (\ref{eq:nls}) and (\ref{NLS2}).  The ground state solution of  $\partial_{xx}Q+Q^5=Q$ reads
\begin{equation}\label{Q def}
Q(x)=\left(\frac{3}{16}\right)^{\frac14}
\frac{2}{\sqrt{\cosh(2x)}}.
\end{equation}

I) With  $\omega_0=1$ and some $T_0>1$, the initial data
\begin{equation*}
\phi(x)=\mathcal{P}_{T_0}R(0)=\frac{1}{\sqrt{T_0}}\fe^{i\frac{x^2}{4T_0}}Q\left(\frac{x}{T_0}\right)
\end{equation*}
leads to a solution of the focusing quintic NLS (\ref{NLS-foc}) that blows up at the chosen $T_0$. Concerning Theorem \ref{thm:main1} and the part (i) of Theorem \ref{thm:main2}, we  take $T_0=1.5$, and we aim to show the  existence of solution for (\ref{eq:nls}) with $u_0=\phi$ or (\ref{NLS2}) with $u_0=\overline{\phi}$ till long times. To do so, we accurately solve the models, and we record the change of the peak value of the solution $|u(t,x)|$ along time, as well as the dynamics of the energy $E_\gamma(t)$ defined in (\ref{energy}).
The results for the DM-NLS model (\ref{eq:nls}) and for the NM-NLS model (\ref{NLS2}) are presented in Figure \ref{fig:1} and Figure \ref{fig:2}, respectively. The profiles of the solutions from different focusing/defocusing time interval are shown in Figure \ref{fig:1a}. From the numerical results, we can observe the following:
\begin{itemize}
    \item For both models, there is no sign of blowup during the dynamics: the solutions remain bounded in $L^\infty$ and so does the energy.
    \item The energy, in both models, is largely increasing in time and  is eventually approaching a  constant. It behaves as a step function, as it is conserved within each half period (each layer) of (\ref{DM}).  The increment is induced at the switching point $t=2n+1$, where the system turns from the focusing layer into the defocusing layer, and the decrement is from the contrary. The  decrement at $t=2n$ is very small, which indicates that at each end of the defocusing layer, the potential energy becomes subtle compared to the kinetic energy.
    \item  The dynamics in the DM-NLS and in the NM-NLS are quite different from each other:   the peak value of the solution is oscillating within a bounded interval in DM-NLS, while it is decreasing in NM-NLS as time gets large. This means that the solution of DM-NLS is bouncing  between a more localized state and a more expanded  state, while the solution of NM-NLS keeps expanding to far field in this test. We can notice this also from their solution  profiles.
   \item The oscillatory profile in the left plot  of Figure \ref{fig:1} indeed meets our definition of dispersion-manageable in (\ref{manageable}) and therefore supports  Conjecture \ref{Conjecture DM} for $(NLS)^\gamma$.
\end{itemize}

The rigorous analysis for all these observations would be quite  interesting but challenging.

\begin{figure}[t!]
\centering
\psfig{figure=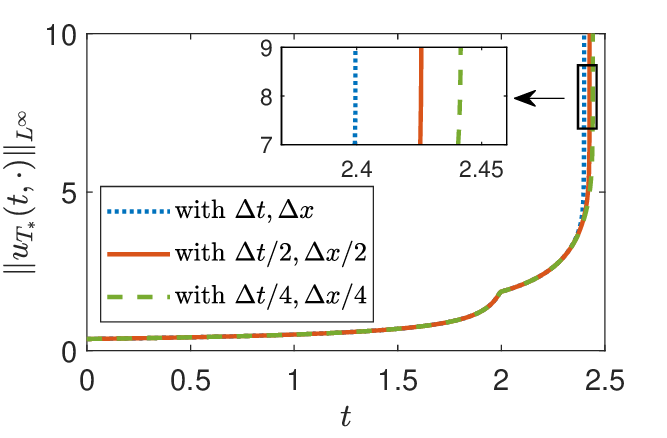,width=0.48\textwidth}
\psfig{figure=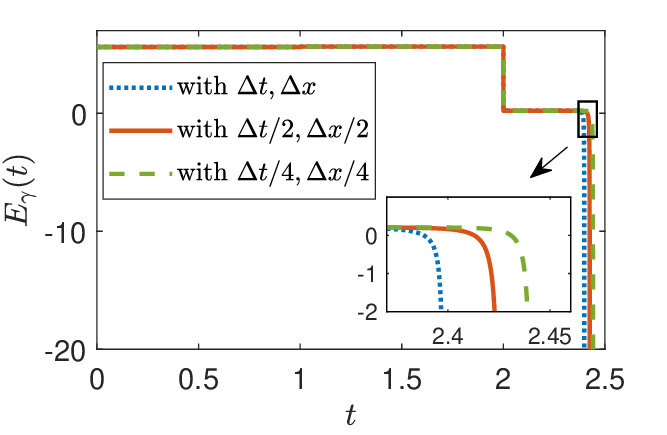,width=0.48\textwidth}
\caption{Blowup in NM-NLS: evolution of the maximum value $\|u_{T_*}(t,\cdot)\|_{L^\infty}$ (left) and the energy $E_\gamma(t)$ (right).
}
\label{fig:thmii}
\end{figure}
\begin{figure}[t!]
\centering
\psfig{figure=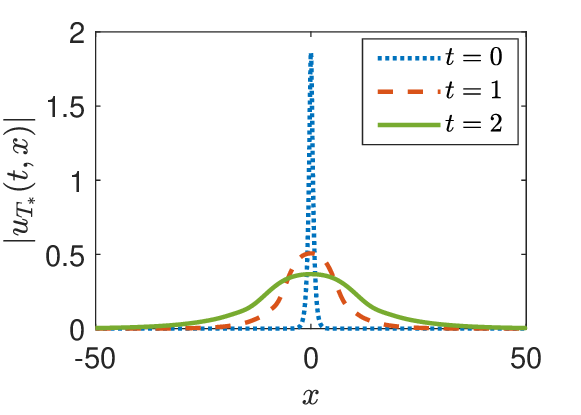,width=0.48\textwidth}
\psfig{figure=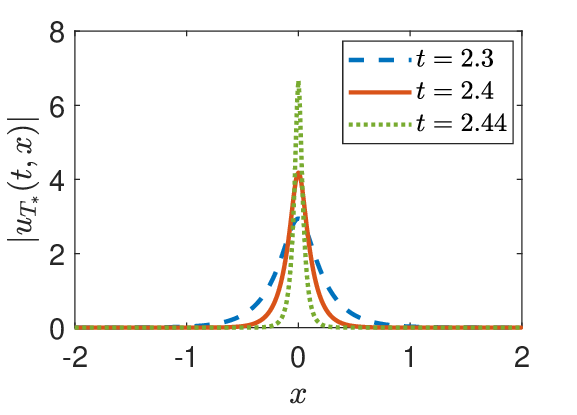,width=0.48\textwidth}
\caption{Blowup in NM-NLS: profiles of the solution $|u_{T_*}(t,x)|$ at different time.
}
\label{fig:thmii2}
\end{figure}

 II) Then, we consider the numerical illustration  for the part (ii) of Theorem \ref{thm:main2}. We take $T_*=2.5$ (i.e., $n=1$) as an example. To find the blowup solution $u_{T_*}$, we follow the construction from the proof in Section \ref{subsection thm2ii}. That is to solve
\begin{equation}
	\left\{ \aligned
	&i\pd_t {\tilde u} +\De \tilde u = \gamma(2-t)| \tilde u|^{4} \tilde u, \quad  0<t\leq2,\nonumber\\
	& \tilde u(0,x) ={\mathcal{P}_{T_*}R(2)},
	\endaligned
	\right.
\end{equation}
and then set $$u_0(x)=\overline{\tilde u(2,x)}.$$
With the numerically obtained $u_0$ from the above, we then solve the NM-NLS (\ref{NLS2}) till
$t=T_*$. Note that for the computation of a blowup solution, one could only get close to the blowup time point but can never reach it \cite{CPC} due to the approximation errors and the limited numerical stability of the numerical scheme. Here to verify the part (ii) of Theorem \ref{thm:main2}, our computations are made by using a second order finite difference in time and Fourier pseudo-spectral in space scheme starting   with the step size $\Delta t=5\times10^{-4},\Delta x=0.046$ and with the refined ones. The changes of $\|u_{T_*}(t,\cdot)\|_{L^\infty}$ and $E_\gamma(t)$ during the computations are shown in Figure \ref{fig:thmii}. The profiles of the solution at different time (under $\Delta t/4,\Delta x/4$) are shown in Figure \ref{fig:thmii2}. From the results, we can see:
\begin{itemize}
    \item The maximum value of the solution increases dramatically when $t$ gets close to $T_*=2.5$. It will eventually exceed the limit of the computer and lead to the crash of the computational program. Such blowup shown in the figure is in fact the numerical blowup which occurs always earlier than the true blowup time, as the solution turns too large and breaks the numerical stability.
    \item The piece-wise energy conservation law is maintained until the numerical blowup happens, where the numerical error becomes large and pollutes the system.
    By refining the mesh size, the numerical blowup will get closer and closer to the real one, but due to the limited computational resource, we could never reach it.
    \item The solution is smooth for $t\in[0,2]$ but becomes sharper and sharper for $t>2$. These observations demonstrate that the constructed $u_0$ can lead to a solution $u_{T_*}$ of (\ref{NLS2}) which blows up at $t=T_{*}$.
\end{itemize}

\begin{figure}[h!]
\centering
\psfig{figure=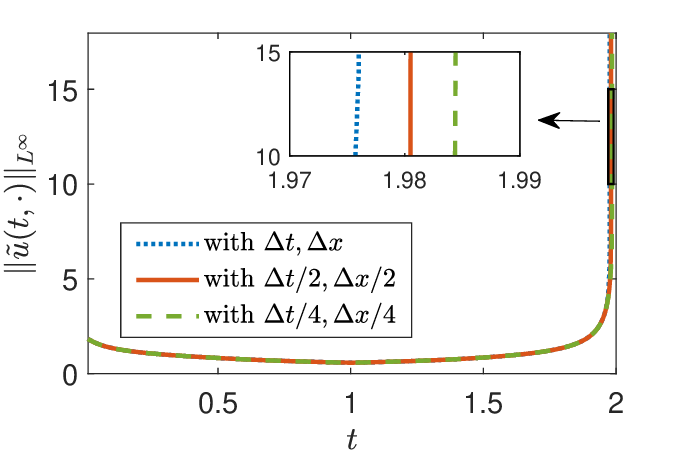,width=4.2cm,height=4cm}
\psfig{figure=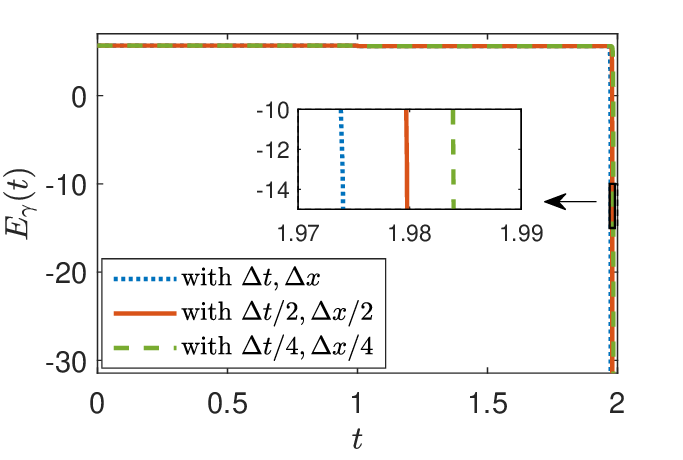,width=4.2cm,height=4cm}
\psfig{figure=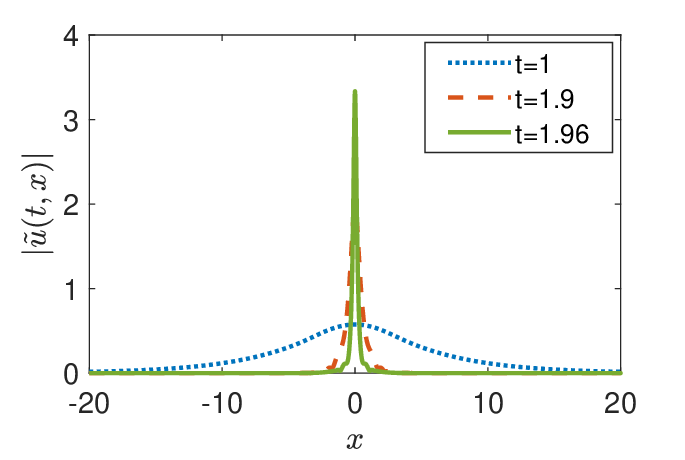,width=4.2cm,height=4cm}
\caption{Results for (\ref{ex eq i}): maximum value (left), energy (middle) and solution profile $|\tilde u(t,x)|$ at different time.
}
\label{fig:exI}
\end{figure}

\subsection{More explorations}
Here we carry out more numerical explorations to address some points that are not covered by our theorems.

I) One may wonder if we can have a blowup solution for DM-NLS as in part (ii) of Theorem \ref{thm:main2}. Here we can demonstrate that at least, a similar   construction as for NM-NLS is not working. More precisely,   there is no $u_0$ for (\ref{eq:nls}) such that $u(x,2n)=\mathcal{P}_{T_*}R(2n)$ with $T_*\in(2n,2n+1)$. To show this numerically, we choose again $d=1$, $n=1$ and $T_*=2.5$, and we solve
\begin{equation}\label{ex eq i}
	\left\{ \aligned
	&i\pd_t {\tilde u} +\gamma(2-t)\De \tilde u = | \tilde u|^{4} \tilde u, \quad  t>0,\\
	& \tilde u(0,x) =\overline{\mathcal{P}_{T_*}R(2)},
	\endaligned
	\right.
\end{equation}
to see if we can reach $t=2$  and obtain $\tilde u(2,x)$. The numerical results are given in Figure \ref{fig:exI}, where $\Delta t=5\times10^{-4}$ and $\Delta x=0.0767$. Clearly, the results suggest that the solution $\tilde u$ of (\ref{ex eq i}) turns to blow up right before $t=2$ and so $\tilde u(2,x)$ does not exist. This  means that one can not find a $u_0(x)$ for the DM-NLS (\ref{eq:nls}) to generate the desired blowup solution $u(2,x)=\mathcal{P}_{T_*}R(2)$. In such sense, we may say that  DM-NLS is more stable than NM-NLS.

\begin{figure}[t!]
\centering
\psfig{figure=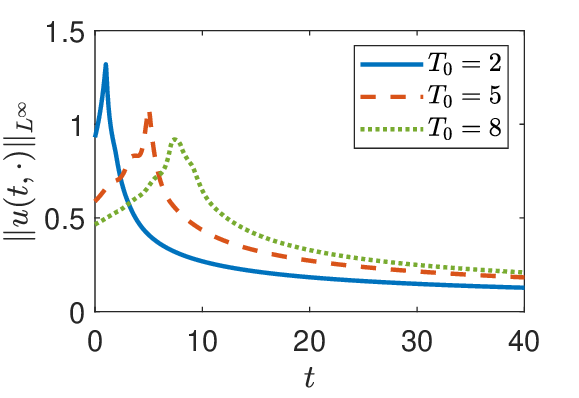,width=0.48\textwidth}
\psfig{figure=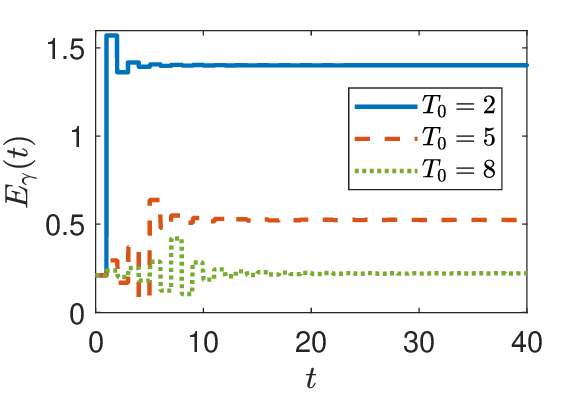,width=0.48\textwidth}
\caption{NM-NLS under  $u_0=\overline{\mathcal{P}_{T_0}R_{\omega_0}(0)}$ with several $T_0\geq2$: the maximum value  (left) and the energy (right) during the dynamics.
}
\label{fig:exii}
\end{figure}
\begin{figure}[t!]
\centering
\psfig{figure=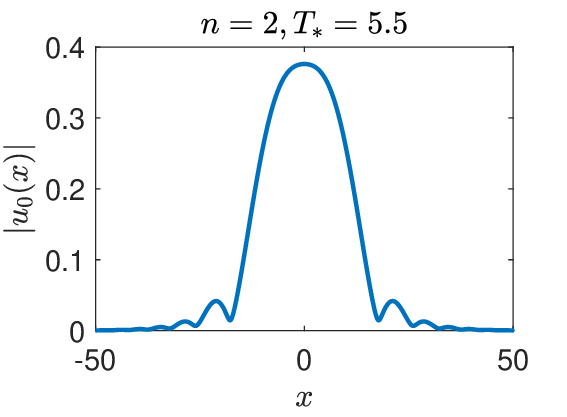,width=0.48\textwidth}\psfig{figure=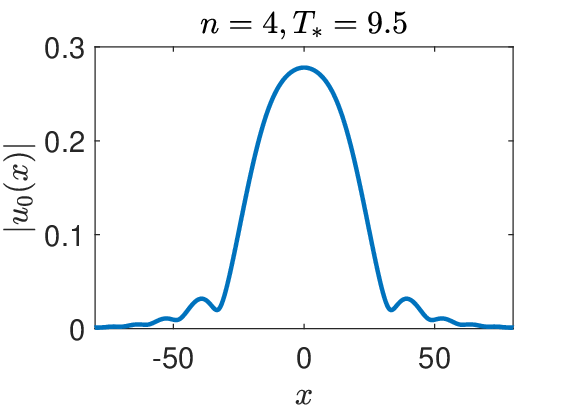,width=0.48\textwidth}\\
\psfig{figure=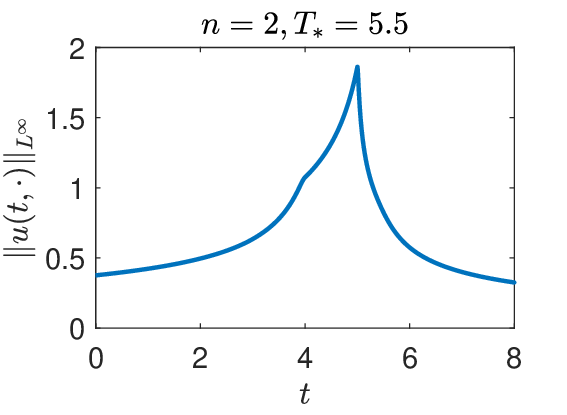,width=0.48\textwidth}
\psfig{figure=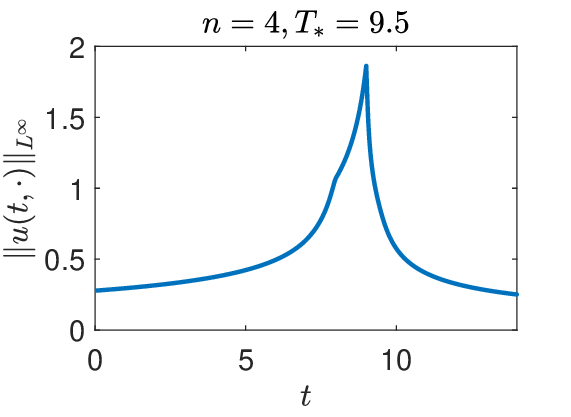,width=0.48\textwidth}
\caption{DM property in NM-NLS:  $u_0$ from (\ref{ex III})  with different $n,T_*$ (1st row) and the maximum value $\|u(t,\cdot)\|_{L^\infty}$ (2nd row) during the dynamics.
}
\label{fig:exiii}
\end{figure}

II) We explore for the case $T_0\geq2$ for the part (i) of Theorem \ref{thm:main2}.  To do so, we take $d=1$ in (\ref{NLS2}) and solve it under  the initial data $u_0=\overline{\mathcal{P}_{T_0}R_{\omega_0}(0)}$ with $T_0=2,5,8$. The dynamics of the maximum value of the solution and the energy are shown in Figure \ref{fig:exii}. From the numerical results, we see the model runs stably and there is no sign of blowups.  Therefore,  the part (i) of Theorem \ref{thm:main2} may hold for all $T_0>1$ which  supports partially the Conjecture \ref{conj 1}, and the rigorous proof will have to be done by other analytical techniques which requires future efforts.

III) We address the Conjecture \ref{Conjecture DM} for the NM-NLS model (\ref{NLS2}). Based on our numerical experience, for most of the initial data (as far as we tested), the intensity of the solution of NM-NLS behaves as the left plot of Figure \ref{fig:1a}, which decreases at large times. However, the part (ii) of Theorem \ref{thm:main2} indeed tells that it can also concentrate at arbitrarily large time. These two facts motive the conjecture. Here we illustrate by showing that for any $n\geq1$, we can find an initial data $u_0(x)$ for (\ref{NLS2}) that leads to a solution satisfying $u(2n,x)=\overline{\mathcal{P}_{T_*}R(2n)}$ with $T_*>2n+1$. To find such $u_0$, we solve
\begin{equation}
	\left\{ \aligned
	&i\pd_t {\tilde u} +\De \tilde u = \gamma(2n-t)| \tilde u|^{4} \tilde u, \quad  0<t\leq2n,\label{ex III}\\
	& \tilde u(0,x) ={\mathcal{P}_{T_*}R(2n)},
	\endaligned
	\right.
\end{equation}
and let $u_0(x)=\tilde u(2n,x)$. Let us take $d=1$ in (\ref{NLS2}), and we choose $n=2,T_*=5.5$ and $n=4,T_*=9.5$. Figure \ref{fig:exiii} shows the obtained $u_0$ and the corresponding dynamics of $\|u(t,\cdot)\|_{L^\infty}$ in the NM-NLS
 (\ref{NLS2}). The results indicate that 1) firstly, the initial data exists and is smooth and localized; 2) the solution of NM-NLS can maintain its intensity up to any prescribed   time interval as long as the $u_0$ is correctly constructed.

 IV) We consider some general cases of blowups in the focusing NLS including a 1D example and a 2D example.

 \begin{figure}[t!]
\centering
\psfig{figure=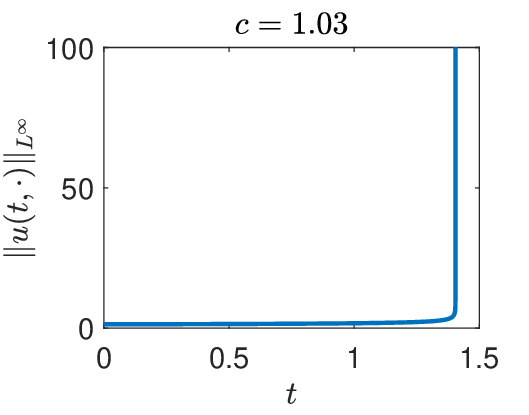,width=0.48\textwidth}
\psfig{figure=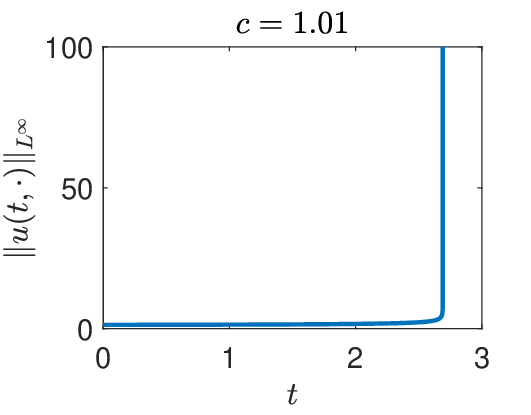,width=0.48\textwidth}
\caption{$\|u(t,\cdot)\|_{L^\infty}$ of focusing NLS (\ref{NLS-foc}) under $u_0=cQ$.
}
\label{fig:exiv foc}
\end{figure}

 \begin{figure}[t!]
\centering
{DM-NLS:}\\
\psfig{figure=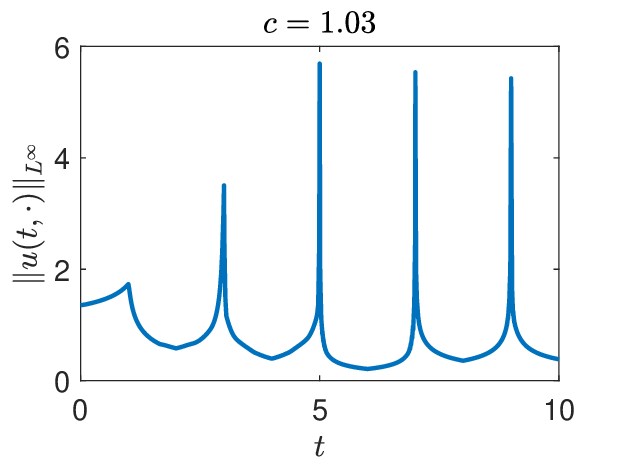,width=0.48\textwidth}
\psfig{figure=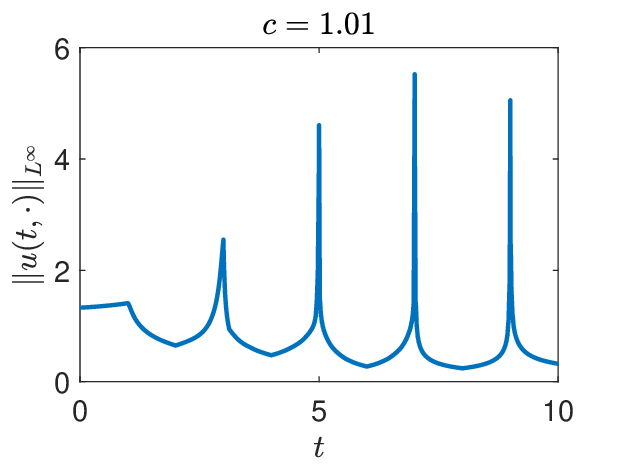,width=0.48\textwidth}\\
{NM-NLS:}\\
\psfig{figure=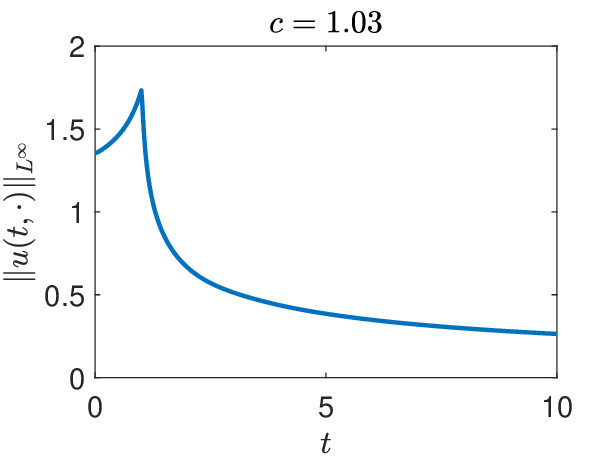,width=0.48\textwidth}
\psfig{figure=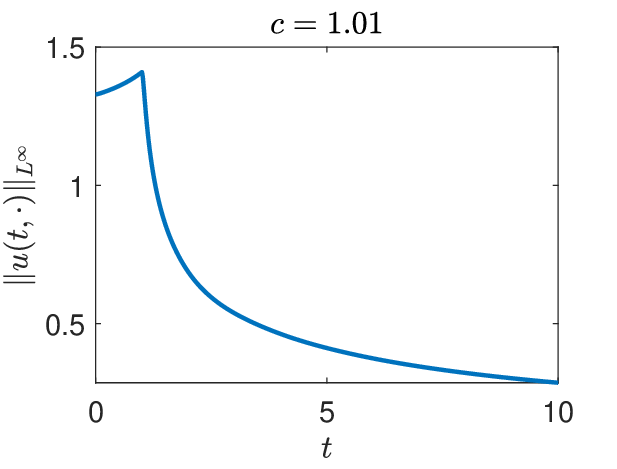,width=0.48\textwidth}
\caption{$\|u(t,\cdot)\|_{L^\infty}$ of DM-NLS and NM-NLS under $u_0=cQ$.
}
\label{fig:exiv DM}
\end{figure}

\begin{figure}[t!]
\centering
\psfig{figure=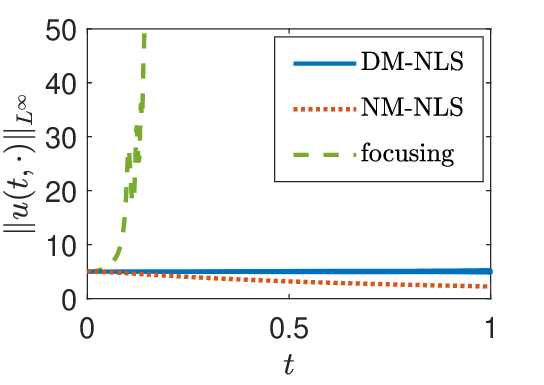,width=0.45\textwidth}
\caption{Evolution of $\|u(t,\cdot)\|_{L^\infty}$ for the 2D example.
}
\label{fig:exIV 2d max}
\end{figure}
\begin{figure}[h!]
\centering
\psfig{figure=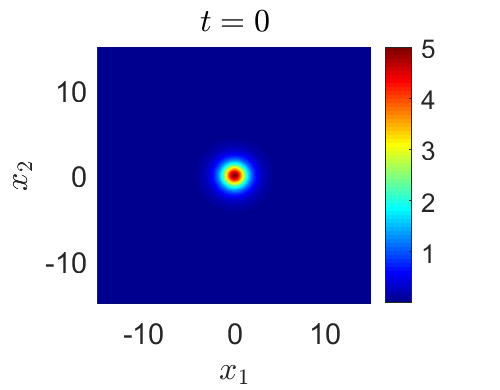,width=0.32\textwidth}
\psfig{figure=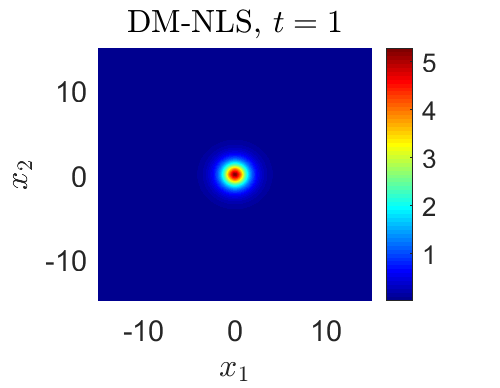,width=0.32\textwidth}
\psfig{figure=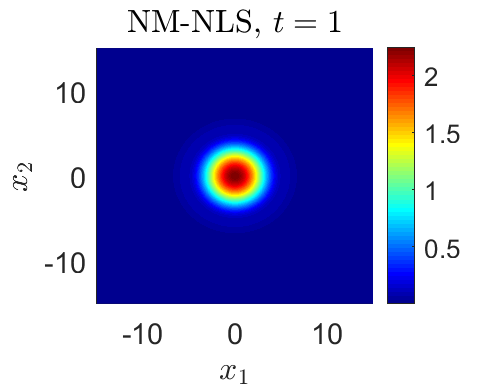,width=0.32\textwidth}
\caption{Solution profile $| u(t,x)|$ for the 2D example.
}
\label{fig:exIV 2d}
\end{figure}

 \begin{itemize}
     \item With $d=1$ and the ground state $Q$ in (\ref{Q def}), we take $u_0(x)=c\,Q(x)$ for some $c>1$ which gives $E_\gamma(0)<0$. By taking $c=1.03$ and $c=1.01$, we can solve the focusing NLS (\ref{NLS-foc}) numerically (with $\Delta t=2.5\times10^{-4},\Delta x=0.0767$) to find that the solution (see Figure \ref{fig:exiv foc}) will blow up at around $t\approx1.4$ and $t\approx2.7$, respectively. Now we input such $u_0$ (it is real-valued) as initial data for both DM-NLS (\ref{eq:nls}) and NM-NLS (\ref{NLS2}) and we simulate the dynamics. The evolution of the maximum value of the solution is shown in Figure \ref{fig:exiv DM}. For  $c=1.03$ or $1.01$, we do not detect any blowups for both models, and their solutions stably pass the original blowup time ($t\approx1.4,\,2.7$) of the focusing NLS. Their dynamics seem to be global.  This means the management (\ref{DM}) can at least delay the   blowup in the original NLS or even suppress it. This   supports the Conjecture \ref{conj 1}.

     \item Finally, we  test a two-dimensional ($ d=2$) example, i.e.,
     $$
i\partial_t u+\gamma(t)\Delta u=|u|^{2}u\ \mbox{ or }\
i\partial_t u+\Delta u=\gamma(t)|u|^{2}u,\quad x=(x_1,x_2)\in\R^2.
$$
With the initial data $u_0(x)=5\mathrm{sech}(|x|/0.86)$, the focusing NLS (\ref{NLS-foc}) forms  (numerically with $\Delta t=2.5\times10^{-5},\Delta x=0.1473$) very quickly a blowup at $t\approx0.14$. See Figure \ref{fig:exIV 2d max}. To stabilize it, we adjust the period for the management $\gamma(t)$ in (\ref{DM}) as
$t_0=0.001,\,t_*=\frac{t_0}{2}$ and keep $\gamma_\pm=1$.  With such $\gamma(t)$, we compute
the dynamics in the above DM-NLS and NM-NLS and the focusing NLS (\ref{NLS-foc}). The evolution of the maximum value of the solution is shown in Figure \ref{fig:exIV 2d max}. The profiles of the solutions are given in Figure \ref{fig:exIV 2d}. From the numerical results, we can see that: 1) the DM-NLS and the NM-NLS have delayed or suppressed the blowup;  2) the intensity of the solution of DM-NLS is well maintained with the adjusted $\gamma(t)$, which partially illustrates the Conjecture \ref{Conjecture DM 2}.
 \end{itemize}

\section{Conclusion}\label{sec:con}
\vskip .5cm

In this paper, we investigated   solutions of two mass-critical NLS equations with the management, including a so-called dispersion-managed NLS and a nonlinearity-managed NLS. The consecutive switching between focusing and defousing
layers makes the dynamics of the two models very different from the classical NLS theory. We prove rigorously that the typical blowup phenomena in standard focusing NLS can  be suppressed to some extend, while  there can also exist a solution that blows up in a  desired focusing layer. Extensive numerical studies are carried out to verify the theoretical results, to explore more stabilization effects of the two models and to illustrate their difference.

\section*{Acknowledgments}
J. Li is supported by NSFC 12126408 and Hunan Provincial Natural Science Foundation
of China 2024JJ5004. C. Ning is supported by Science and Technology Program of
Guangzhou No. 2024A04J4027.
X. Zhao is supported by Natural Science Foundation of Hubei Province No. 2019CFA007 and NSFC 42450275, 11901440.
We would like to acknowledge the private communications with J.-C. Saut and J. Murphy.

\section*{Reference}

\bibliographystyle{model1-num-names}

\begin{thebibliography}{00}

\bibitem{NMadd0}
{\sc F.Kh. Abdullaev, J.G. Caputo, R.A. Kraenkel, B.A. Malomed},
Controlling collapse in Bose-Einstein condensates by temporal modulation of the scattering length,
Phys. Rev. A 67 (2003) p. 013605.

\bibitem{2d-dm}
{\sc F.Kh. Abdullaev, B.B. Baizakov, M. Salerno}, Stable two-dimensional dispersion-managed soliton,
Phys. Rev. E 68 (2003) p. 066605.

\bibitem{Saut}
{\sc P. Antonelli, J.-C. Saut, Ch. Sparber}, Well-Posedness and averaging of NLS with time-periodic dispersion management, Adv. Differential Equations 18 (2013) pp. 49-68.


\bibitem{Berge}
{\sc L. Berg\'{e}, V.K. Mezentsev, J.J. Rasmussen, P.L. Christiansen, Y.B. Gaididei}, Self-guiding light in layered nonlinear media, Opt. Lett. 25 (2000) pp. 1037-1039.

\bibitem{Biswas-Book}
{\sc A. Biswas, D. Milovic, M. Edwards}, Mathematical Theory of Dispersion-Managed Optical Solitons,
Nonlinear Physical Science. Springer, Berlin, 2010.

\bibitem{Malomed-kdv}
{\sc S. Clarke, B.A. Malomed, R. Grimshaw},
``Dispersion management" for solitons in a
Korteweg-de Vries system, Chaos 12 (2002) pp. 8-15.

\bibitem{NM-theory1}
{\sc T. Cazenave, M. Scialom}, A Schr\"odinger equation with time-oscillating nonlinearity,
Rev. Mat. Univ. Complut. Madrid 23 (2010) pp. 321-339.

\bibitem{NM-theory2}
{\sc L. Di Menza, O. Goubet}, Stabilizing blow up solutions to nonlinear Schr\"odinger equations,
Commun. Pure Appl. Anal. 16 (2017) pp. 1059-1082.

\bibitem{NM stabilize}
{\sc R. Driben, B.A. Malomed, M. Gutin, U. Mahlab},
Implementation of nonlinearity management for Gaussian pulses in a fiber-optic link by means of second-harmonic-generating modules, Opt. Commun. 218 (2003) pp. 93-104.

\bibitem{GVD1}
{\sc I.R. Gabitov, S.K. Turitsyn}, Averaged pulse dynamics in a cascaded transmission system with passive dispersion compensation, Opt. Lett. 21 (1996) p. 327.



\bibitem{Zhao}
{\sc Y. He, X. Zhao}, Numerical methods for some nonlinear Schr\"odinger equations in soliton management,
J. Sci. Comput. 95 (2023) p. 61.

\bibitem{HeZhao}
{\sc Y. He, X. Zhao}, Numerical integrators for dispersion-managed KdV equation,
Commun. Comput. Phys. 31 (2022) pp. 1180-1214.

\bibitem{CPC}
{\sc X. Hong, Q. Wei, X. Zhao},
Comparison of different discontinuous Galerkin methods based on various reformulations for gKdV equation: soliton dynamics and blowup, Comput. Phys. Comm. 300 (2024)  p. 109180.

\bibitem{Feshbach0}
{\sc S. Inouye, M.R. Andrews, J. Stenger, H.J. Miesner, D.M. Stamper-Kurn, W. Ketterle}, Observation of Feshbach resonances in a Bose-Einstein condensate, Nature 392 (1998) p. 151.


\bibitem{dm-numer1}
{\sc T. Jahnke, M. Mikl},
Adiabatic midpoint rule for the dispersion-managed nonlinear Schr\"odinger equation, Numer. Math. 138 (2018) pp. 975-1009.

\bibitem{dm-numer2}
{\sc T. Jahnke, M. Mikl},
Adiabatic exponential midpoint rule for the dispersion-managed nonlinear
Schr\"odinger equation, IMA J. Numer. Anal. 39 (2019) pp. 1818–1859.

\bibitem{Feshbach1}
{\sc P.G. Kevrekidis, G. Theocharis, D.J. Frantzeskakis, B.A. Malomed}, Feshbach resonance management for Bose-Einstein condensates, Phys. Rev. Lett. 90 (2003) p. 230401.


\bibitem{Malomed-book}
{\sc B.A. Malomed}, Soliton Management in Periodic Systems, Springer, New York, 2006.


\bibitem{Merle-BU-1993}
{\sc F. Merle}, Determination of blow-up solutions with minimal mass for nonlinear Schr\"odinger equations with critical power, Duke Math. J. 69 (1993) pp. 427-454.


\bibitem{Murphy}
{\sc J. Murphy, T. Hoose},
Well-posedness and blowup for the dispersion-managed nonlinear Schr\"odinger equation, Proc. Amer. Math. Soc. 151 (2023) pp.  2489-2502.

\bibitem{VMPG1}
{\sc G.D. Montesinos, V.M. P\'{e}rez-Garc\'{\i}a, P.J. Torres},
Stabilization of solitons of the multidimensional nonlinear
Schr\"odinger equation: matter-wave breathers, Phys. D. 191 (2004) pp. 193-210.


\bibitem{NM1}
{\sc C. Par\'{e}, A. Villeneuve, P.-A. Belange, N.J. Doran}, Compensating for dispersion and the nonlinear Kerr effect without phase conjugation, Opt. Lett. 21 (1996) p. 459.

\bibitem{Saito}
{\sc H. Saito, M. Ueda},
Dynamically stabilized bright solitons in a two-dimensional Bose-Einstein condensate, Phys. Rev. Lett. 90 (2003) p.  040403.

\bibitem{GVD2}
{\sc M. Suzuki, I. Morita, N. Edagawa, S. Yamamoto, H. Taga, S. Akiba}, Reduction of Gordon-Haus timing jitter by periodic dispersion compensation in soliton transmission,  Electron. Lett. 31 (1995) pp. 20-27.

\bibitem{Taobook}
{\sc T. Tao}, Nonlinear Dispersive Equations. Local and Global Analysis. Amer. Math. Soc. Providence 2006.

\bibitem{NM2}
{\sc I. Towers, B.A. Malomed}, Stable $(2+1)$-dimensional solitons in a layered
medium with sign-alternating Kerr nonlinearity, J. Opt. Soc. Am. B 19 (2002) p. 537.



\end{thebibliography}

\end{document}